\documentclass[a4paper, leqno, intlimits]{amsart}

\usepackage{graphicx}
\usepackage{amsmath}
\usepackage{amssymb}
\usepackage{amsxtra}
\usepackage{mathrsfs} 

\usepackage[applemac]{inputenc}

\newcommand{\psh}{\mathcal{PSH}}
\newcommand{\spsh}{\mathcal{SPSH}}

\newcommand{\C}{\mathcal{C}}

\newcommand{\zbar}{\bar{z}}

\newcommand{\dc}[1]{d^c\mkern .8mu #1}
\newcommand{\ddc}[1]{{dd^c}\mkern .8mu #1}
\newcommand{\ma}[2][n]{(\ddc{#2})^{#1}}    

\newcommand{\bd}{\partial}

\newcommand{\CC}{\mathbf{C}}
\newcommand{\RR}{\mathbf{R}}
\newcommand{\ZZ}{\mathbf{Z}}

\newcommand{\Proof}{\noindent\textit{Proof.}\hspace{0.6em}}
\newcommand{\QED}{\qed \vspace{\baselineskip}}

\newcommand{\ie}{i.{\kern0.11111em}e.\@}
\newcommand{\Ie}{i.{\kern0.11111em}e.\@}
\newcommand{\Eg}{E.{\kern0.11111em}g.\@}
\newcommand{\eg}{E.{\kern0.11111em}g.\@}
\newcommand{\cf}{cf.\@}

\theoremstyle{plain}
\newtheorem{theorem}{Theorem}
\newcounter{Lemma}
\newtheorem{lemma}[Lemma]{Lemma}

\author{Jonas Wiklund} 
\title[Monge-Amp\`ere measure at the boundary]{Monge-Amp\`ere measure at
the boundary of some domains with corners}
\date{March 31, 2006}

\begin{document}

\begin{abstract}
    Let $\mu^z$ be the measure obtained by sweeping out the
    Monge-Amp\`ere measure of the pluricomplex Green function with pole
    at $z. $ We prove that $\mu^z$ vanish on Levi flat parts
    of the boundary for 1) every relatively compact analytic
    polyhedron in complex space, 2) product domains of hyperconvex sets
    in Stein manifolds.
\end{abstract}

\maketitle

\section{Introduction}

Let us denote the plurisubharmonic (henceforth abbreviated as psh)
functions on a complex manifold \( X \) by \( \psh(X), \) strictly
psh functions by \( \spsh(X) \) and non-positive psh functions by
\( \psh^-(X).  \)

For notations of the complex Monge-Ampère measure at the boundary of
an hyperconvex manifold we refer to \cite{Demailly:85, Demailly:87}. Let
us just recall some basic definitions.

Suppose \( \varphi \) is a psh weight on a hyperconvex manifold \(
\Omega \), \ie\ \( \varphi \in \psh^-(\Omega), \) the pseudo balls \(
B_r = \{ z \::\: \varphi(z) < r\} \) are relatively compact in \(
\Omega \) for some \( r \in \RR, \) and \( e^\varphi \) is a continuos
function.  We denote the pseudo spheres by \( S_r = \{ z \::\:
\varphi(z) = r\}.  \)

Denote the characteristic function of a set \( A \) by \(
\mathbf{1}_A. \) Suppose \( X \) is Stein, for any psh weight \( u \)
on \( X \) Demailly defined a ``swept out'' measure of \( \ma{u} \) on \(
S_r \) by
\begin{equation}\label{eq:definition-of-mu}
   \mu_{u,r} = \mathbf{1}_{X \setminus B_r} \ma{u} -
   \ma{(\max(u,r))}.
\end{equation}
For smooth plurisubharmonic \( u \) the formula
\( 
   \mu_{u,r} = \ma[n-1]{u} \wedge \dc{u} \big|_{S_r}
\) holds (\cf\ \cite{Demailly:85}).
Using continuity under monotonic sequences  the measure \(
\mu_u \) for a psh weight \( u \) on a hyperconvex domain \( \Omega \)
is defined
by \[ \mu_u = \lim_{r \nearrow 0} \mu_{u,r}.  \]

In \( \CC^n \) we may define the \textit{Lelong number} of a psh function \( u \)
at the point \( w \) as 
\begin{equation} 
    \nu(u,x) = \lim_{r \to 0} \frac{M(u,r,w)}{\log{r}},
\label{eq:lelongdef2}
\end{equation}
where \( M(u,r,w) = \sup\{u(z)\::\: |z-w| < r \}.  \) Since the Lelong
number of a psh function is well known to be independent of a
biholomorphic change of variables the definition readily carry over to
any complex manifold.

The \textit{pluricomplex Green function} with \textit{pole at} \( w \)
was introduced by Klimek \cite{Klimek:85} and
Zahariuta~\cite{Zahariuta:84}.  For any connected relatively open
subset \( \Omega \) of \( X \) it is defined as
\[ 
   g_{\Omega}(z,w) = \sup\{u(z)\;;\; u\in \psh^-(\Omega), \; \nu(u,w) 
   \geq 1 \}.
\]

Finally, fix a pole \( w \) in \( \Omega,  \) and define a measure on the
boundary by
\[ 
   \mu^w = \lim_{r\to 0}
   \mu_{g_{\Omega}(z,w),r}.
\]




In the original paper by Demailly following geometric classification
theorem was obtained.

\begin{theorem} 
    \emph{(Demailly \cite{Demailly:87}).  }\/ For any domain \( \Omega
    \) with defining psh function \( \rho \in \C^2(\overline{\Omega}) \) the
    measure \( \mu^{w} \) is supported on the strict pseudoconvex
    points of the boundary, for every \( w \in \Omega.  \)
\end{theorem}

In this note we give two examples that this theorem may be generalized
to some domains with corners.

\section{Hyperconvex domains with corners}

Let us start with a smoothing Lemma, similar to Lemma 3.2 in
\cite{Guan:02}. 
\begin{lemma}\label{lem:Guan-smoothing}
    Let \( U \) be an open domain in a complex space \( X. \) Suppose
    \( u_i \in \psh\cap \C^k(U), i=1,2, \) for some \( k \in
    \ZZ_{\geq 0} \) and let \( \chi_r: \RR \to \RR \) be
    a smooth convex function such that \( \chi_r(x) = |x| > r.  \)
    Define \( M_r(u_1,u_2) = (\chi_r(u_1-u_2) + u_1 +
    u_2)/2 \).  Then \( M_r(u_1,u_2) \in \psh(U) \) and \( M_r(u_1,u_2) 
    \geq 
    \max(u_1,u_2) \)  with \( M_r(u_1,u_2) 
    =
    \max(u_1,u_2) \) except on a small neighbourhood of the
    ar\^ete \( \{|u_1|=|u_2|\}.  \) Furthermore  \( M_r \in \C^k.  \)
\end{lemma}

Guan proved this Lemma for twice differentiable functions by direct
computation of the Levi form for \( M_r.  \) In \cite{Cegrell:01b} a
much simplified proof of this Lemma was given, where the Lemma was
stated for all psh functions \( u_i \) regardless of continuity; and it
is not clear from Cegrell's proof whether \( M_r \) is upper
semicontinuous.  For the convenience of the reader we repeat the
argument here.

\Proof  Since \( 2\max(u_1,u_2) =
|u_1-u_2| + u_1+u_2 \) the latter part is clear.  For the first part
note that
\[ 
   \chi_r(x) = \sup\{kx+l, |k| \leq 1, \text{ and } kt+l \leq
   \chi_r(t), t \in \RR \},
\]
thus
\[ 
   \chi_r(u_1-u_2) + u_1 + u_2 = \sup\{(1+k)u_1 + (1-k)u_2 + l \},
\]
and \( (M_r)^* \) is psh, but since \( u_i \) are continuous we
have that \( (M_r)^* = M_r. \) 
\QED

\begin{theorem}
    Suppose \( P \) is a relatively compact analytic polyhedron, and
    let \( w \in P \), then \( \mu^w \) is supported in a subset of
    the strictly pseudoconvex points of \( \bd P. \)
\end{theorem}

\Proof
Let \( P \) be an analytic polyhedra defined by
\[ 
   P = P(f_1,\ldots,f_N) = \{ z \in W \::\: |f_i| < 1, i = 1,\ldots,N \}
\]
where \( W \Subset G \subset \CC^n, \) and \( f_i \in
\mathscr{O}(G,\CC)  \) is not identically zero. P is clearly
hyperconvex.

Define \( \tilde{\varphi} := \max_k \{\log|f_i|\}, \) then \(
\tilde{\varphi} \) is a psh weight on \( P, \) furthermore \(
\exp(\tilde{\varphi}) \in \C(\bar{P}). \) 

From the defining Equation~\eqref{eq:definition-of-mu} for the swept out
Monge-Ampère measure it is clear that \( 
g_P(z,w) \) and \( \max\{g_P(z,w),-1\} \) have the same boundary
measure on the pseudo spheres \( \{ z\::\:g_P(z,w) = t \} \) for \( t >
-1, \) thus they have the same boundary measure on \( \partial P. \)

By the continuity and maximality of \( g_P \) there is a positive
constant \( c \) such that \( \max\{g_P(z,w),-1\} \geq c
\tilde{\varphi}(z), \forall z \in \overline{\Omega}.  \) Let \(
\varphi = c\tilde{\varphi}. \)
Thus
by the comparison principle for the boundary measures (\cf\ Theorem
3.4 \cite{Demailly:87}) we have 
\begin{equation}\label{eq:comp}
    \mu^w \leq \mu_{\varphi}. 
\end{equation}

Now a straightforward calculation gives that
\begin{multline}\label{eq:calculus}
    (\ddc \log|f_i|)^{n-1} \wedge \dc{\log|f_i|} =
    \frac{1}{(2|f_i|^2)^n} (\ddc |f_i|^2)^{n-1} \wedge \dc{|f_i|^2} = 
    \\
    \frac{1}{(2|f_i|^2)^n} (\ddc{|f_i|^2-1})^{n-1} \wedge
    \dc{(|f_i|^2)}.
\end{multline}
Take  a point \( x \) in a Levi flat part of the boundary of \( P \) then \(
|f_i|^2-1 \) is a local defining function for \( P \) around \( x, \) 
thus the last expression in Equation~\eqref{eq:calculus} vanish.

Take any point \( x \in \bd P, \) away from the corners of \( P. \)
Take \( R > 0 \) fixed and set \( u_R =
M_R(\log|f_1|,\ldots,\log|f_N|), \) where we have extended \( M_R \)
to \( N \) variables in the natural way.

Since \( u_R \) is smooth
on \( P \) we have in \( U_x \cap P, \) where \( U_x \) is some
neighbourhood of \( x, \)
\[ 
   \mu_{u_R}|_{S_t} = (\ddc{u})^{n-1} \wedge \dc{u} =
   \frac{1}{(2|f_i|^2)^n} (\ddc{|f_i|^2-t})^{n-1} \wedge
       \dc{(|f_i|^2)} = 0.
\]
By continuity under decreasing sequences we have that \( \lim_{R \to 0}
\mu_{u_R}|_{S_t} = \mu_u |_{S_t} = 0 \) and then letting \( t \nearrow 0 \) we
may conclude that \( \mu_{\varphi} = 0 \) on Levi flat parts of the
boundary of \( P. \)
By Equation~\eqref{eq:comp} \( \mu^w \leq \mu_{\varphi} = 0 \)
 on Levi flat part of \( \bd P. \)
\QED

\begin{theorem}\label{thm:cross}
    Suppose \( \Omega_i \) is a hyperconvex domain in a 
    Stein manifold \( X_i, \) of dimension \( n_i, \) for \( i =1,2,
    \) and let \( X = X_1 \times X_2, \) and \( \Omega = \Omega_{1}
    \times \Omega_{2} \subset X \).  Take \( w \in \Omega, \) then \(
    \mu^{w} \) vanish except on the corner \( \bd \Omega_1 \times \bd
    \Omega_2.  \)
\end{theorem}

\Proof
Fix the pole \( w = (w_1, w_2) \in \Omega_1 \times \Omega_2, \)
and let \( u_1\) and \( u_2 \) be the pluricomplex Green function of
\( \Omega_1 \) and \( \Omega_2 \) with poles at \( w_1 \) and \( w_2
\) respectively.

Let \( u = \max\{u_1,u_2\}.  \) Then \( u \) is a psh weight by the
continuity result in \cite{Demailly:87}, and therefore \( \Omega \) is
a hyperconvex set in \( X. \) Clearly \( u = 0 \) on \( \bd \Omega.
\) Directly from Equation~\eqref{eq:lelongdef2} we have \( \nu(u,w) =
\min\{\nu(u_1,w_1),\nu(u_2,w_2)\} = 1.  \) Since \( \ma[n_i]{u_i} = 0
\) away from \( w_i, \) for \( i=1,2, \) we have by Proposition 3.4 in
\cite{Zeriahi:91} that \( \ma{u} = 0 \) outside the pole \( w.  \) By
\cite{Demailly:87, Zeriahi:96} \( u \) is in fact the pluricomplex
Green function for \( \Omega, \) with pole at \( w.  \)

Now, for any twice differentiable function \( v \) depending only on
\( n < n_1+n_2 = N \) variables we have either that \( n < N-1 \)
which implies that \( \ma[N-1]{v} = 0, \) or we have \( n = N-1.  \)
If \( n = N-1 \) we get, in local coordinates, \( \ma[N-1]{v} = f(z)\,
dz_1 \wedge d\zbar_1 \wedge \ldots \wedge dz_{N-1} \wedge
d\zbar_{N-1}, \) for some function \( f.  \) Thus
\begin{multline}\label{eq:N-1}
   \ma[N-1]{v} \wedge \dc {v} = \\ 
   = f(z)\, dz_1 \wedge d\zbar_1 \ldots\wedge
   dz_{N-1} \wedge d\zbar_{N-1} \wedge \Big(\sum_{k=1}^{N-1}
   \frac{\partial v}{\partial z_k}dz_k - \frac{\partial v}{\partial 
   \zbar_k}d\zbar_k\Big) = 0.
\end{multline}

Fix \( r<0.  \) Take \( \rho_i \in\spsh\cap \C(X_i), \) \( i=1,2 \)
then there is a decreasing sequence of smooth psh functions \(
\{u_i^j\}_{j=1}^{\infty} \) on \( \Omega_i \) such that \(
\rho_i/(j+1)+u_i < u_i^j < \rho_i/j + u_i, \) since the \( u_i \)'s
are at least continuous.  Thus on a neighbourhood of \( S_r \) we have
a decreasing sequence of psh functions \( u^j = \max\{u_1^j,u_2^j\},
\) smooth away from a neighbourhood of the arête \( \{|u_1|=|u_2|\}.
\) But by Equation~\eqref{eq:N-1} we have \( \ma[N-1]{u^j} \wedge \dc
u^j \equiv 0, \) away from the arête.

By continuity of the currents under decreasing limit we have \(
\mu_{u,r} = 0.  \) Thus, away from the arête \( \{ z \in \Omega \::\:
|u_1(z)| = |u_2(z)| \} \) we have \( \mu_w |_{S_r} = 0.  \) Letting \(
r \to 0 \) gives the result.  \QED

Let us just add the remark that the proof of Theorem~\ref{thm:cross}
could be made a lot cleaner by using the definition of the boundary
measure from an unpublished manuscript of Cegrell~\cite{Cegrell:PC05}.
Alas, the boundary measure is only defined on hyperconvex set in \(
\CC^n \) in that paper.

\end{document}